\documentclass[12pt,oneside,english]{amsart}
\usepackage[latin1]{inputenc}
\usepackage{amssymb}

\makeatletter

\usepackage{babel}

\makeatother
\begin{document}
\baselineskip=17pt

\title{New digit results and several problems}

\author{Vladimir Shevelev}
\address{Departments of Mathematics \\Ben-Gurion University of the
 Negev\\Beer-Sheva 84105, Israel. e-mail:shevelev@bgu.ac.il}


\begin{abstract}
We give some new relations for Newman digit sums respectively
different modulos and put some problems. In particular, for the odd
prime  modulos we put an important conjecture.
\end{abstract}

\maketitle
\section{Introduction}
As in \cite{7} put for $q>1$

\begin{equation}\label{1}
n=\sum^v_{k=0}a_kq^k,\;\;\; 0\leq a_k< q,\;\;\;
\sigma_q(n)=\sum^v_{k=0}a_k.
\end{equation}

Denote for $x\in\mathbb{N},\;\;\;l\in[0,m-1]$

\begin{equation}\label{2}
S_{m,l,q}(x)=\sum_{0\leq n < x: n\equiv l(\mod
{m})}(-1)^{\sigma_g(n)}
\end{equation}

In the case $q=2$ we write
$S_{m,l,2}=S_{m,l},\;\;\sigma_2(n)=\sigma(n)$.

We call (\ref{2}) a generalized Newman sum.

In \cite{7} we gave a quite another proof of the Coquet's estimates
for $S_{3,0}(x)$ and a fast algorithm for its calculation. Professor
J.-P.Allouche kindly informed me about a misprint in Coquet's
theorem: for odd $x$

$$
\eta=\eta(x)=(-1)^{\sigma(3x-1)}
$$

(but $(-1)^{\sigma(3x-3)}$ as in \cite{2}; sf. \cite{1}, pp.98-99)

An important role in our proof belongs to the formula: for an even
$n$

\begin{equation}\label{6}
S_{3,0}(4n)=3 S_{3,0}(n).
\end{equation}

The method of proof (\ref{6}) in \cite{7} allows to obtain several
new relations for  some Newman digit sums and to formulate a very
important conjecture.

\newpage

\section{Some new digit relations}

We use the following simple relations for
$S_{m,l}(x),\;x\in\mathbb{N}.$

If $m$ is odd then

$$
S_{m,l}(2x)=\begin{cases} S_{2m,l}(2x)+S_{2m,l+m-1}(2x),\;\; if \;\;
l\;\; is \; even\\-S_{2m,l-1}(2x)+ S_{2m, l+m}(2x), \; if \;\; l\;
\; is \; odd\end{cases}=
$$

\begin{equation}\label{7}
=\begin{cases} S_{m,\frac l 2}(x) - S_{m, \frac{l+m-1}{2}}(x), \; if
\;\; l \;\; is\;\; even\\
-S_{m,\frac{l-1}{2}}(x)+S_{m,\frac{l+m}{2}}(x), \;\; if \;\; l \;\;
is \; odd.\end{cases}
\end{equation}

If $m$ is even then

\begin{equation}\label{8}
S_{m,l}(2x)=\begin{cases} S_{\frac m 2,\frac l 2}(x),
\; if \;\; l \;\; is\;\; even\\
-S_{\frac m 2,\frac{l-1}{2}}(x), \; if \;\; l \;\; is \;
odd\end{cases}.
\end{equation}

Note that (\ref{8}) reduces the calculations to the case of an odd
$m$. Hence, for an odd $m$ we should solve the system (\ref{7}) to
get an equation for e.g. $S_{m,0}(x)$ only.

The calculations by this method are rather long and sometimes
complicated. Nevertheless, we obtained the following relations for
$x,y\in\mathbb{N}$, the first of which was obtained in \cite{7}
(here as in \cite{7} $S_{m,0}([y, y+z))=S_{m,0}(y+z)-S_{m,0}(y)$):

\begin{equation}\label{9}
S_{3,0}([8x,\; 8y))= 3 S_{3,0}([2x,\; 2y)),
\end{equation}

\begin{equation}\label{10}
S_{5,0}([32x,\; 32y))= 5 S_{5,0}([2x,\; 2y)),
\end{equation}

\begin{equation}\label{11}
S_{7,0}([128x,\; 128y))= -7 S_{7,0}([2x,\; 2y)),
\end{equation}

$$
S_{9,0}([512x,512y))= 3 S_{9,0}([128x,128y))+
$$

\begin{equation}\label{12}
+3S_{9,0}([8x,8y))- 9 S_{9,0}([2x,\; 2y)).
\end{equation}

Besides, by similar way we obtained the following relation for
$S_{5,0,4}([x,y)):$ if $x$ is divisible by $32$ then
\begin{equation}\label{13}
S_{5,0,4}([256x,256y))= 10 S_{5,0,4}([16x,16y))- 5 S_{5,0,4}([x,y)).
\end{equation}

\newpage

Using (\ref{9})-(\ref{13}) as in \cite{7} it could be proved that

\begin{equation}\label{14}
|3S_{5,0}(n)|=O(n^{\frac{\ln 5}{\ln 16}})= O(n^{0,58048})
\end{equation}

\begin{equation}\label{15}
|S_{7,0}(n)|=O(n^{0.46789\ldots}),
\end{equation}

\begin{equation}\label{16}
|S_{9,0}(n)|=O(n^{0.79248\ldots}) (as \;\; for \;\; S_{3,0}(n)),
\end{equation}

\begin{equation}\label{17}
|S_{5,0,4}(n)|=O(n^{0.81092\ldots}).
\end{equation}

\section{Some conjectures and problems}

1) To find a method (probably, a variant of the method of generating
functions) for an automatic obtaining of relations of type
(\ref{9})-(\ref{13}). To find a general digit equation of this type
(at least, for the base 2).

2) According to (\ref{9})-(\ref{11}) we have  in particular that

\begin{equation}\label{100}
S_{3,0}(2^3 )= 3,\;\;S_{5,0}(2^5 )= 5,\;\;S_{7,0}(2^7 )= -7.
\end{equation}

Denote

$$
a_n=S_{n,0}(2^n).
$$

By the further direct calculations for the prime values of $n$ we
obtained a very astonishing sequence:

$$
a_3=3,\;a_5=5,\;a_7=-7,\;a_{11}=11,\;a_{13}=13,
$$
\begin{equation}\label{101}
a_{17}=697,\;a_{19}=19,\;a_{23}=-23,\;a_{29}=29, \ldots
\end{equation}

It was very difficult for us to believe that $a_{17}=697$!

It this connection recall a remarkable result of M.Drmota and
M.Skalba \cite{3}: the only primes $p\leq 1000$ satisfying
$S_{p,0}(n) >0$ (at least, for sufficiently large $n$) are
$3,5,17,43,257,683$.

Therefore it is natural to conjecture that for primes \slshape
different \upshape from $17,43,257,683,\ldots$ we have

\newpage

\begin{equation}\label{102}
S_{p,0}(2^p)=\pm p.
\end{equation}

Note that, (\ref{102}) satisfies also for $3$ and $5$  because of
the numbers $2^3$ and $2^5$ are small.

Besides we conjecture that \slshape always \upshape $p |
S_{p,0}(2^p)$.

3). In the connection with the results (\ref{10}),(\ref{11}) it is
interesting to find the sharp estimates in  these cases similar to
\cite{2} and \cite{7}.

4)In our opinion, it is very interesting to find a generalization of
(\ref{13}) for $S_{2k+1,0,2k}(x)$ and get the sharp estimates.

We conjecture that not only $S_{2k+1,0,2k}(x)> 0,\;\;k\geq 1$, but
also the Newman-like phenomena becomes more and more strong with the
enlargement of $k$ . Moreover, if

$$
S_{2k+1,0,2k}(x)=O(x^{\lambda_k})
$$

then we conjecture that $\lim_{k\rightarrow\infty}\lambda_k=1$.

5) We conjecture that, if $d|m$ then the characteristic polynomial
which corresponds to the relation of considered type for
$S_{m,0,q}(x)$ is divisible by one for $S_{d,0,q}(x)$.

6) We conjecture that if $(m,3)=1$ then for any $k\in(1,\frac m 3)$,

$$
|S_{m,0}(x)|=o(|S_{3k,0}(x)|).
$$

Note that, if 6) is true then it could be proved our Conjecture 2
\cite{5} which until now has only a heuristic justification
\cite{6}.

\bfseries Remark. \mdseries   Put
\begin{equation}\label{20}
G^{(i)}_{m,0}(x) = \sum_{0\leq n < x, n \equiv 0(mod \; m),
\sigma(n)\equiv i (mod \; 2)} 1,\;\;\;\; i=0,1.
\end{equation}

It is a special case of the Gelfond digit sum.  It is evident that

$$
G^{(o)}_{m,0}(x)+G^{(1)}_{m,0}(x)=\sum 1=\lfloor\frac x m\rfloor+1,
$$

$$
G^{(o)}_{m,0}(x)-G^{(1)}_{m,0}(x)=S_{m,0}(x).
$$
\newpage

By the Gelfond theorem

\begin{equation}\label{21}
G^{(i)}_{m,0}(x)=\frac {x}{2m}+O(x^{\frac{\ln{3}}{\ln{4}}}).
\end{equation}

Thus, estimates (\ref{10}), (\ref{11})  make more precise the
remainder term in (\ref{21}) in the cases of  $m=5\;and\;m=7$.

\section{Acknowledgements}

The author is grateful to Professor D.Berend for the placing at his
disposal of the paper \cite{3} and to Professor J.-P.Allouche for
information of a Coquet's misprint and indication the reference
\cite{1}.

\end{document}